\RequirePackage{fix-cm}
\documentclass[twocolumn]{svjour3}          % twocolumn
\smartqed  % flush right qed marks, e.g. at end of proof
\usepackage{graphicx}
\usepackage{amssymb,amsmath}
\usepackage{pgfplots}
\usepackage{booktabs}
\usepackage{appendix}
\usepackage{xcolor}
\usepackage{multirow,tabularx}
\usepackage[position = top, justification=raggedright,singlelinecheck=false]{subfig}
\usepgfplotslibrary{colormaps}
\usetikzlibrary{pgfplots.colormaps}
\pgfplotsset{compat=1.16}

\usepackage[square,sort,comma,numbers]{natbib}
\usepackage{url}
\begin{document}
\sloppy

\title{Robustness of nonlinear parameter identification in presence of process noise using control-based continuation%\thanks{Grants or other notes
}


\author{Sandor Beregi         \and
        David A. W. Barton \and
        Djamel Rezgui \and
        Simon A. Neild %etc.
}


\institute{S. Beregi \at
              University of Bristol, Faculty of Engineering \\
              \email{sandor.beregi@bristol.ac.uk}           %  \\
           \and
           D. A. W. Barton \at
                         University of Bristol, Faculty of Engineering\\
                         \email{david.barton@bristol.ac.uk}           \\
           D. Rezgui \at
                         University of Bristol, Faculty of Engineering \\
                         \email{djamel.rezgui@bristol.ac.uk}           \\
           S. A. Neild \at
                         University of Bristol, Faculty of Engineering \\
                         \email{simon.neild@bristol.ac.uk}           %  \\
}

\date{Received: date / Accepted: date}

\maketitle

\begin{abstract}
In this study, we consider the experimentally-obtained, periodically-forced response of a nonlinear structure in the presence of process noise. Control-based continuation is used to measure both the stable and unstable periodic solutions while different levels of noise are injected into the system. Using this data, the robustness of the control-based continuation algorithm and its ability to capture the noise-free system response is assessed by identifying the parameters of an associated Duffing-like model. We demonstrate that control-based continuation extracts system information more robustly, in the presence of a high level of noise, than open-loop parameter sweeps and so is a valuable tool for investigating nonlinear structures.
\end{abstract}

\section{Introduction}

Studying physical structures experimentally can be a challenge if the measurements are polluted with a significant amount of noise. Often information is lost and it becomes difficult to resolve the fine details of the experiment's behaviour. In this paper, we propose control-based continuation as an approach for robust parameter identification in noisy nonlinear systems. By tracking the steady-state solutions, we analyse a noise-contaminated experimental rig featuring nonlinear characteristics and assess the performance of the control-based continuation algorithm benchmarked against open-loop parameter sweeps.

Investigating nonlinear systems is a long-established field of dynamical analyses \cite{GuckHolm}. Nonlinear dynamical models have been used effectively to explain phenomena in many engineering applications \cite{nodycon_proc} as well as in many other areas such as, amongst others, biological \cite{biology} or economical \cite{economics} systems.

From a practical point of view, an essential part of analysing nonlinear phenomena is building models which represent the relevant features of the system. This is often achieved using reduced-order mathematical models, which requires the identification of the model parameters.

In the literature, there are numerous studies advocating a range of different methods for parameter identification. In engineering applications, a common solution is to measure the restoring force of the system and directly identify the system parameters as shown by Masri et al \cite{Masri2004} and Cammarano et al \cite{Cammarano2011}. While this approach can be effective in characterising the system, the required measurements can be difficult to carry out in some cases. Therefore, many papers focus on methods that identify model parameters based on the observed dynamical behaviour. Kerschen et al. \cite{Kerschen2006} and Noel and Kerschen \cite{Noel2017} give a thorough review covering a range of different techniques of nonlinear system identification. A number of these methods, such as nonlinear auto-regressive moving average with exogenous inputs (NARMAX) modelling \cite{Leontaritis1, Leontaritis2} or the use of the Hilbert transform are based on the analysis of time series \cite{Fledman}. Goharoodi \cite{Goharoodi2018} and Marchesiello \cite{Marchesiello2008} perform nonlinear subspace identification by a time-domain study of the system response for a given excitation. A similar approach is followed by Narayanan \cite{Narayanan} using multi-harmonic excitation. Noel and Kerschen carry out a similar identification of nonlinear subspaces in the frequency domain \cite{Noel2013} reducing the computational burden of the method \cite{Noel2014}.

Another possible strategy is to extend the theory of (linear) modal analysis to nonlinear systems by considering nonlinear normal modes \cite{Shaw, Vakakis} and measuring the nonlinear frequency response diagram of the structure \cite{carrella_FRF, Noel2016}.

{ In some cases, when the nonlinearity is weak, a system may have a quasi-linear frequency response at small forcing amplitudes with a unique and stable periodic solution corresponding to every forcing frequency. This may allow the use of linear techniques in the analysis e.g. to filter out noise \cite{Link, Karaagacli}. Nevertheless, this study addresses the problem of tracing steady-state solutions in parameter ranges when significant nonlinear behaviour (bistability, hysteresis) can be observed.}

In practice, the { nonlinear} frequency-response of a system can be obtained in several ways for example by applying random excitation or impulses to the system. In \cite{Peter1, Peter2} Peter et al. obtain the backbone curve and frequency response by controlling the phase-shift between the excitation and the response. In our study, we focus on the commonly-used approach where the frequency-response is extracted by performing parameter-sweeps. However, nonlinear systems may have specific (so-called bistable) parameter domains where two (or even more) stable steady-state solutions exist for the same set of parameters. The boundary of the domains of attraction of these stable solutions is referred to as a separatrix which can be often characterised by another, unstable, steady-state solution. In an experiment, one can observe this as a sensitivity of the steady-state behaviour to the initial conditions. For example, one may find that a certain level of perturbation is tolerable and so the dynamics stay within the same basin of attraction, while a larger perturbation may lead the system to diverge from its originally observed steady-state behaviour. This phenomenon also means that by standard parameter-sweeps, only stable solutions can be captured. Thus, part of the bifurcation diagram may remain hidden. The presence of process noise can also contribute to the loss of data by potentially driving the system away from a stable solution with a small domain of attraction. Unlike measurement noise, which can be largely eliminated thanks to the current sensor and data acquisition technology, process noise is an inherent feature of the observed phenomenon, i.e. certain effects that are too complex for deterministic modelling are best captured by a stochastic process. In these cases it is an interesting challenge to identify the underlying deterministic dynamics as its features may be disguised by the random perturbations.

%An additional challenge of these algorithms, compared to linear parameter-identification, is that while in nonlinear systems a complicated structure of stable and unstable steady-state solutions may exist, in experiments, only stable solutions can be found directly while usually limited information can be acquired about the unstable ones.

Continuation methods trace a family of solutions in a nonlinear system by applying small parameter changes to follow the solution branch. In case of mathematical models, several analytical or numerical techniques, based on bifurcation analysis, are available to track steady-state solutions irrespective of their stability \cite{Kuznetsov}. Numerical  continuation methods, in particular, are constructed around solving a nonlinear zero problem. For example, periodic solutions are commonly calculated using the method of collocation to discretise the solution over the time-period while the solution branch is continued with the pseudo-arclength method with respect to a system parameter, referred to as bifurcation parameter \cite{Seydel}.

Control-based continuation \cite{Sieber2008} is a method which incorporates the techniques of numerical continuation and bifurcation analysis to trace solutions of physical and numerical experiments where the governing equations are not explicitly available. Thus, in principle, it is capable of capturing both stable and unstable steady-state solutions. Just as is the case for numerical bifurcation analysis, periodic solutions have to be discretised and composed into an equivalent zero problem, e.g. by the coefficients of their truncated Fourier series. These coefficients can be used to design a control which is stabilising and non invasive; that is, the steady-state solutions of the controlled system are the same as the steady-state response of the open-loop system.

In our study, we apply the method of control-based continuation to conduct the experimental bifurcation analysis of a forced nonlinear oscillator subjected to noise. The acquired data is then used to identify parameters for a model of the structure based on the assumption that it is a one-degree-of-freedom Duffing-like oscillator, which we use as a model of the experimental rig. { While one could use the measured time profiles for the identification, here, we stick to the S-shaped amplitude response curves to assess the robustness of control-based continuation to noise directly on the bifurcation diagrams, which are frequently in the focus of studies on nonlinear systems.}

By physically polluting the system with different levels of noise, we assess the robustness of control-based continuation derived data to random perturbations. This is an important factor from the point of view of the practical effectiveness of the method as, in general, a heavy noise-load may result in losing relevant information from the system. The robustness of control-based continuation is tested and compared against parameter sweeps without control. Thus, we assess if control-based continuation is capable to capture the response more accurately and if it is capable to reveal details from the system which would otherwise remain undetectable, providing a more robust basis for model building and parameter identification.

\section{Motivation}

\begin{figure}
  \centering
    \includegraphics{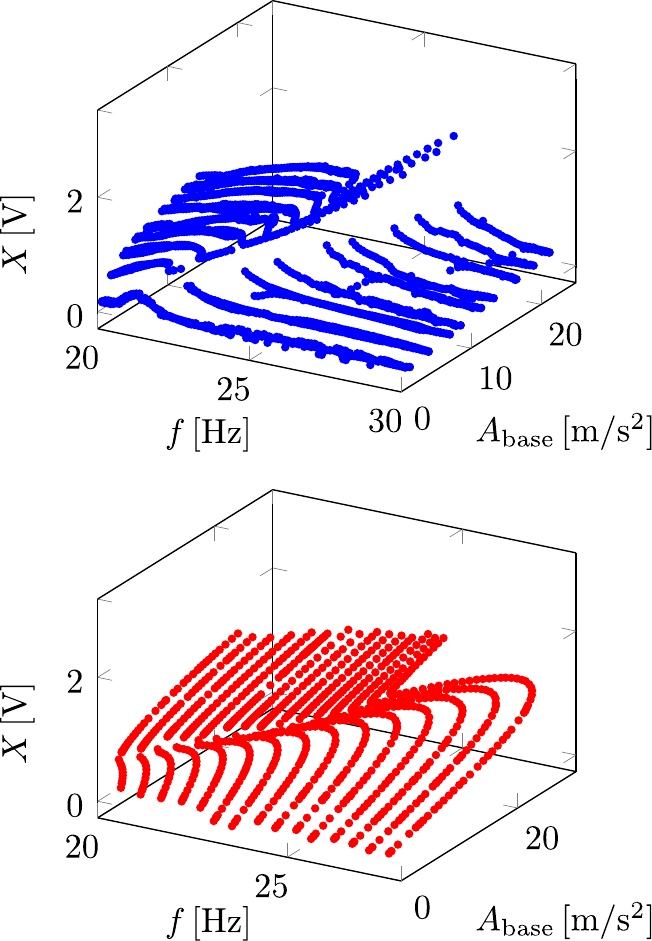}
\caption{Frequency sweeps (top panel), and control-based continuation (bottom panel).}\label{fig:sw_cbc_surf}
\end{figure}

In case of nonlinear structures, it is a common approach to use the frequency response to characterise the system, see e.g. \cite{Agarwal2018, Tang2016}. This response is often measured in a forced or base-excitation setting, by performing two frequency sweeps, one with increasing and another with decreasing forcing frequency. Thus, in principle, every stable solutions can be captured even in bistable frequency domains. With the help of frequency-sweeps, performed at different forcing amplitudes, one can trace the stable part of the response surface in the forcing frequency -- forcing amplitude parameter plane.

An example is shown in the left panel of Fig. \ref{fig:sw_cbc_surf} with the response surface of the vibration amplitudes $X$ for the nonlinear oscillator, described in Section 3. It can be observed that, while the frequency-sweeps were performed at constant shaker voltage amplitude due to the internal dynamics of the shaker, the shaker acceleration amplitude $A_{\rm base}$ varies slightly during each sweep. We traced the same surface by performing control-based continuation at constant frequencies (see the right panel). This method is not only capable of capturing the unstable solutions but the feedback control on the periodic solution also results in a smoother response surface.

If the experiment is subject to process noise, it can be challenging to trace steady-state solutions in parameter ranges where the system undergoes bistable behaviour. This bistable behaviour is typically characterised by the presence of three solution branches within the same parameter range, where two solution branches are connected by an unstable branch through saddle-node bifurcations. If parameter sweeps are performed, in addition to the lack of information regarding the unstable solutions, one may erroneously observe that the vibration amplitude jumps earlier than the saddle-node bifurcations are located in the solution branch of the underlying noise-free system. Ultimately, this leads to a loss of information not only about the unstable branch but also parts of the stable solution branches. Moreover, it may not be possible to identify a clear fold point in the captured solution manifold as the noise may cause the system to jump repeatedly between the domains of attraction of two steady-state solutions \cite{Barkley}.

\begin{figure*}[]
    \centering
    \includegraphics{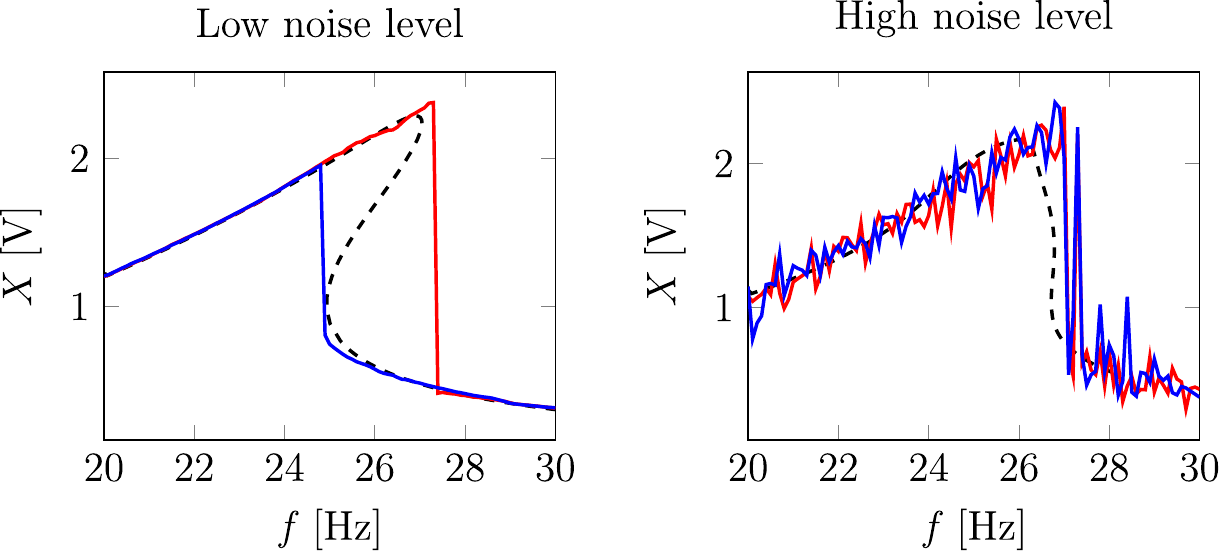}
    \caption{Comparison of frequency sweeps (red and blue curves) in case of low and high levels of process noise. A seventh-order Duffing-like oscillator best fit is indicated by the dashed black curves.}
    \label{tab:freq_sweep_noise}
\end{figure*}

This phenomenon is illustrated by the example in Fig. \ref{tab:freq_sweep_noise} where two frequency sweeps, carried out on the nonlinear oscillator, are compared. In both cases, the oscillator was modelled as a seventh-order Duffing-like oscillator. { The model was fitted to measurement data using the analytical formulae for the frequency response.} The response of the fitted model is indicated by the black curves. It can be seen that in the low noise case, in spite of not having information about the unstable solutions, the fitted mathematical model matches with the measured vibration amplitudes reasonably well. However, with higher levels of noise, the bistable domain is not properly traced and the fitted model provides an entirely different response to that using the low-noise data. { Note that the parameter-fitting method used in this motivating example is the same as that used in the discussion below. The reader may refer to Section 5 and the Appendix for details.}

 %Such behaviour was reported and investigated in connection with many engineering applications, such as in case of the oscillation of towed wheels, flutter instability of aeroplane wings or the chatter instability of machine-tools.

\section{The experimental rig}

\begin{figure}
\begin{center}
\includegraphics[width=80 mm]{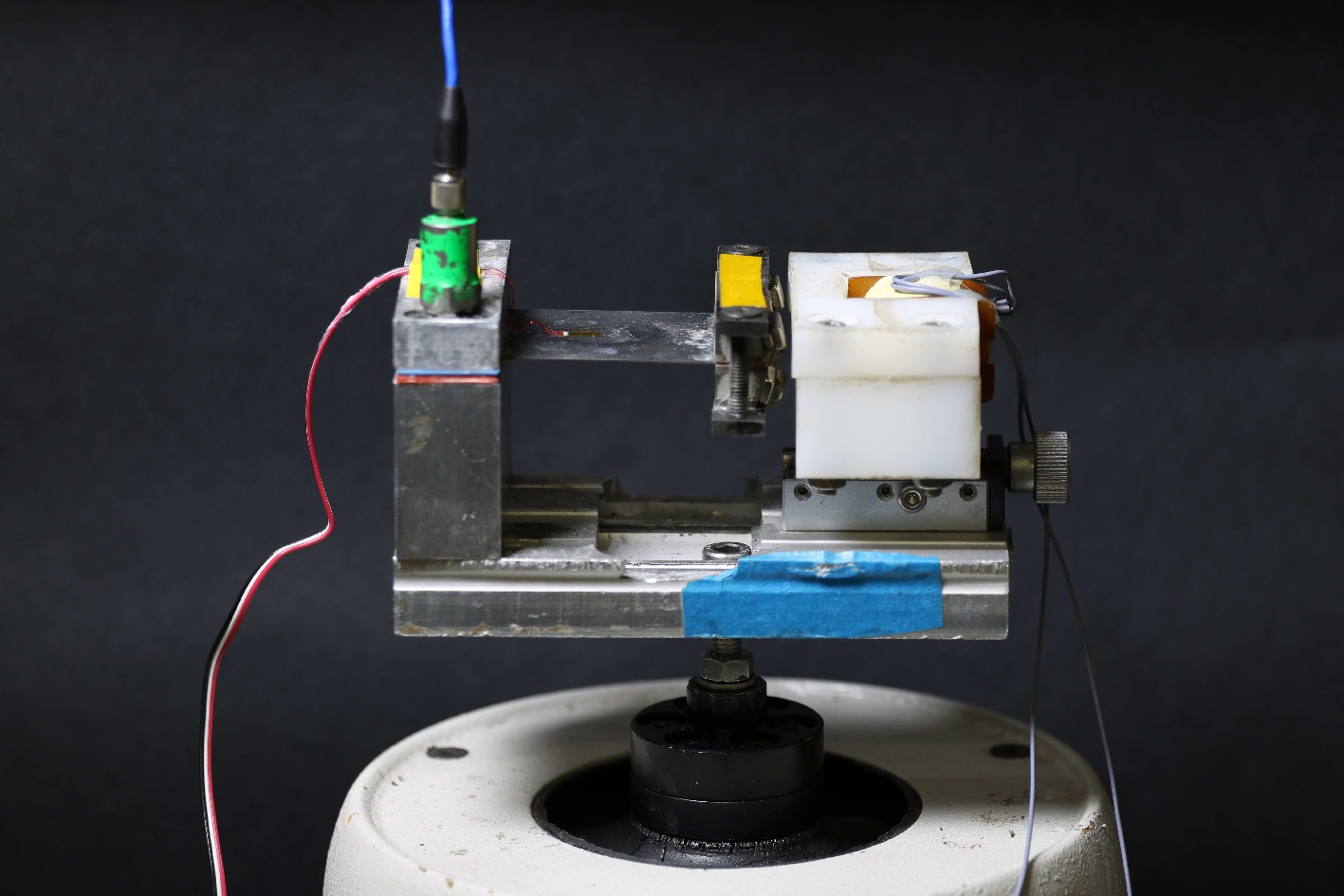}
\caption{The experimental setup: a nonlinear oscillator mounted on the shaker.}\label{fig:meas_photo}

\end{center}
\end{figure}

\begin{figure}
\begin{center}

\includegraphics{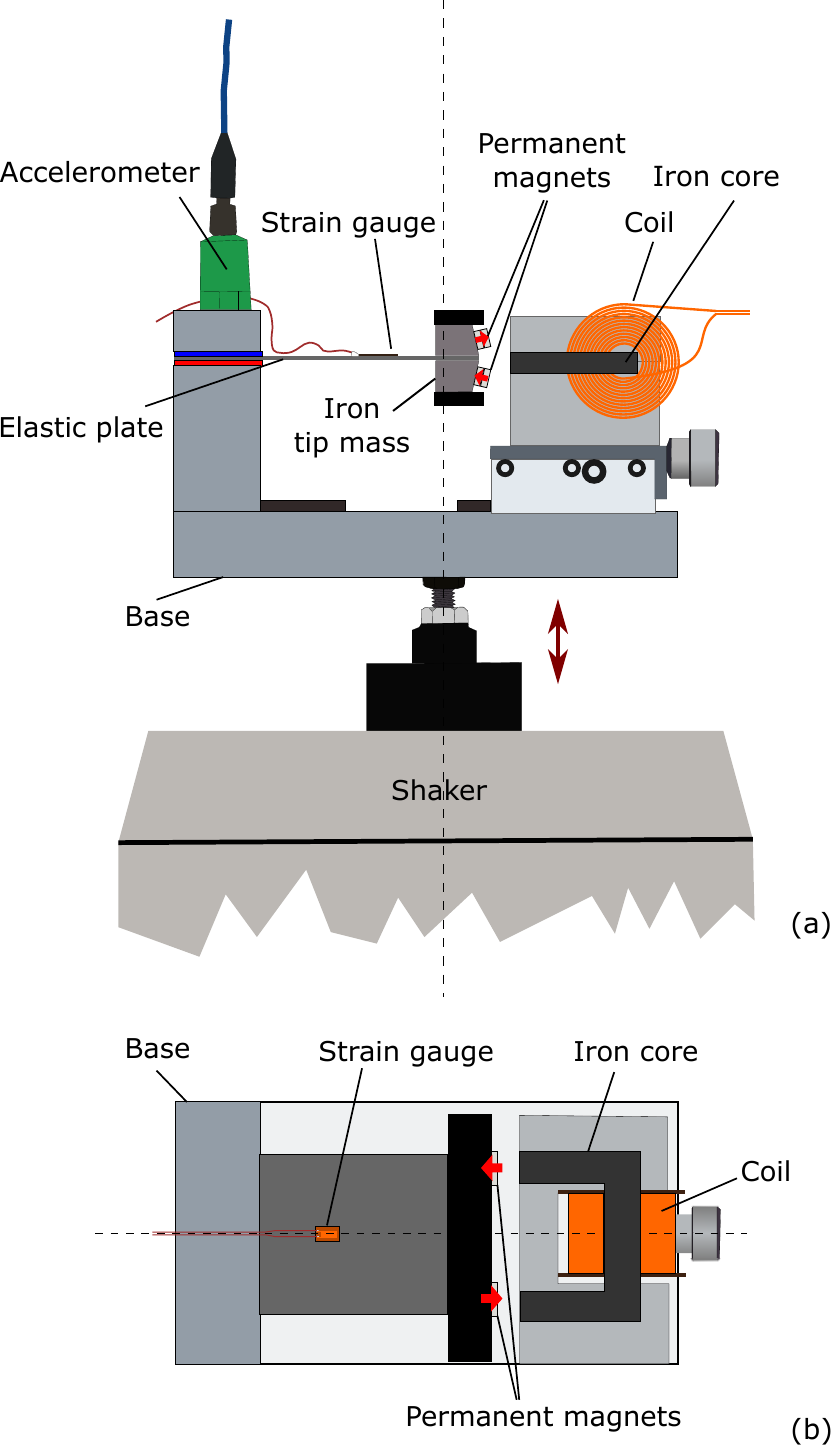}
\caption{The experimental rig: schematics of the nonlinear oscillator on a shaker. Panel (a): side view, panel (b) top view. For visibility of the coil and its core, their housing is shown in the background only.}\label{fig:exp_rig}

\end{center}
\end{figure}

\begin{figure*}
\begin{center}
\includegraphics{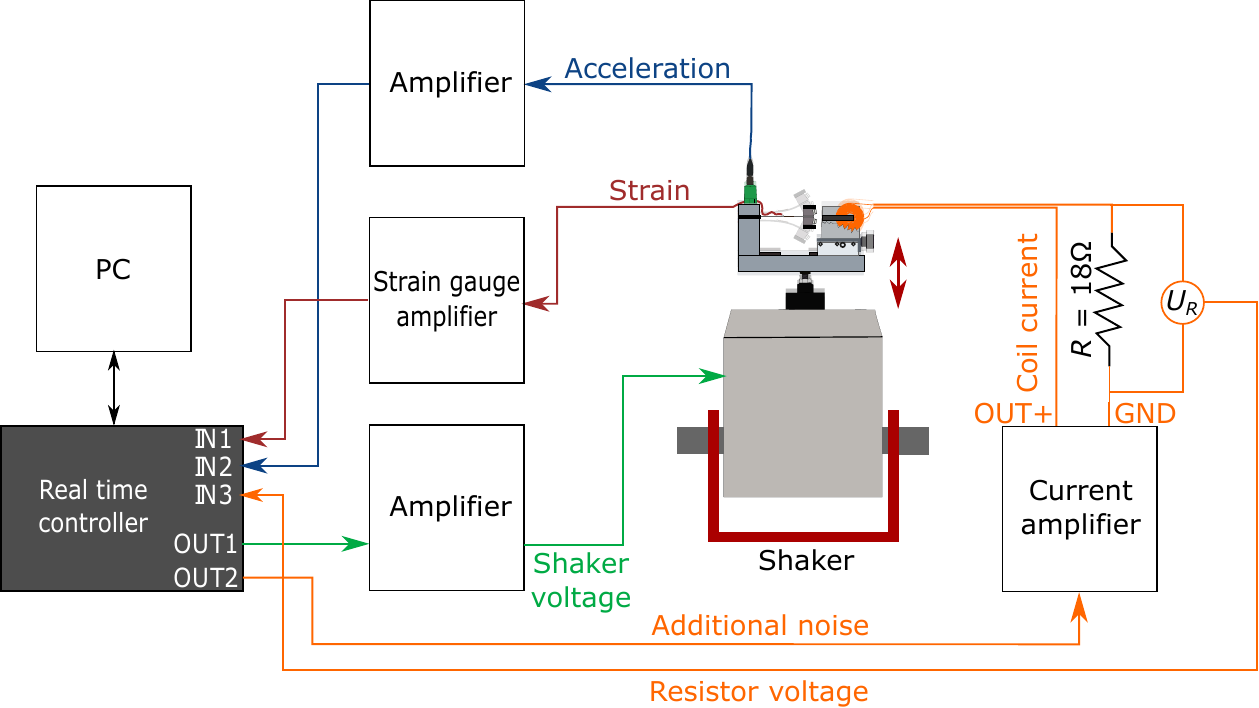}
\caption{The experimental rig: schematics of the experimental setup for control and data acquisition.}\label{fig:exp_rig_overall}
\end{center}
\end{figure*}

In our study, we investigate the properties of the forced nonlinear oscillator shown in Fig. \ref{fig:meas_photo}. This device is designed to show nonlinear behaviour to make it suitable to test the capabilities of experimental algorithms \cite{Barton2011, Barton2012}, while in former studies, it was also used as an energy harvester \cite{Cammarano2011}. The structure of the device is shown in Fig. \ref{fig:exp_rig}. The nonlinear oscillator is formed from a thin steel plate, which is clamped to the base as a cantilever beam. This clamp includes plastic plates supporting the steel plate to add damping to the contact. At the other end, two iron masses, incorporating four permanent magnets, are attached to the tip of the plate. The orientation of the poles of these magnets is indicated by the red arrows in Fig. \ref{fig:exp_rig}. The resultant magnetic field interacts with a stator (see the white block in Fig. \ref{fig:meas_photo}), which consists of an electromagnetic coil with an iron core in an insulated housing.

The rigid base of the oscillator and the stator is mounted on a shaker, providing forcing to the system. As the plate vibrates, the magnets at the free end will move relative to the coil changing the magnetic flux. The resulting combination of structural, inertial, and magnetic forces results in a nonlinear restoring force. According to Faraday's law of electromagnetic induction, the changing magnetic flux also generates an electromotive force, which can be measured as a changing voltage across the coil. Similarly, a current flowing through the coil results in a magnetic flux and so generates a magnetic force, which allows the excitation of the system through the coil. In our experiments, we used the coil to pollute the experiment with additional noise.

Figure \ref{fig:exp_rig_overall} shows a schematic of the experimental set up including data acquisition and control. The data acquisition and demand generation is carried out by a real-time controller box built upon a BeagleBone Black board \cite{DBrtc}. Three input and two output voltage signals are handled by the real-time controller. The demand voltage OUT1 is applied to the shaker using an amplifier while the demand OUT2 is used to contaminate the experiment with noise. For this purpose, we use a current amplifier to regulate the current flowing through the coil according to the noisy signal generated by the controller box.

The vibration of the elastic plate is measured using a strain-gauge, which is connected to an amplifier in a quarter-bridge configuration. The voltage, provided by the strain-gauge amplifier, is measured at the acquisition channel IN1. Channel IN2 is connected to an accelerometer that is mounted on the top of the base of the nonlinear oscillator, while IN3 records voltage that is generated across a resistor by the current flowing through the coil (see the circuit in Fig. \ref{fig:exp_rig_overall}).

\section{Control-based continuation}

We employ the technique of control-based continuation to determine the response of the nonlinear oscillator to periodic forcing. Let us consider a one-degree-of-freedom forced nonlinear system as a model of the experimental rig in the form of
\begin{equation}\label{eq:1dof_eom}
\ddot x + g(\dot x, x) = F(t),
\end{equation}
where the state variable $x$ is the input voltage from the strain-gauge (IN1), the dot refers to differentiation with respect to time, the function $g$ contains all the nonlinearities, while $F(t)$ corresponds to the forcing provided by the shaker. Despite the fact that this model does not include the electromagnetic effects in the system explicitly, we found that the one-degree-of-freedom model can characterise the response to periodic forcing with good accuracy. Note, that, as demonstrated in \cite{Cammarano_thesis}, this is not the case for non-periodic excitation.

The forcing $F(t)$ is composed by two parts: a periodic component and an additional control term $F_{\rm ctrl}$
\begin{equation}
F(t) = A \cos(\omega t) + B \sin(\omega t) + F_{\rm ctrl},
\end{equation}
where $\omega$ is the angular frequency of the forcing whereas $A$ and $B$ are constant coefficients.
To capture the nonlinear response of the open-loop (uncontrolled) system, the control has to fulfil two conditions: it has to be stabilising and non-invasive, i.e. the steady-state response of the controlled system has to be equal to the steady-state response of the open-loop system. This is ensured by the proportional-derivative control law given in the form
\begin{equation}\label{eq:control}
F_{\rm ctrl} = k_{\rm p}(x^\ast - x) + k_{\rm d}(\dot x^\ast - \dot x),
\end{equation}
where $x^\ast$ is the control target while $k_{\rm p}$ and $k_{\rm d}$ are the proportional and derivative control gains. { Since the derivative $\dot x$ is not acquired directly in our experiment, its value at a given time $t_i$ is obtained using a backward difference scheme
\begin{equation}
\dot x (t_i) \approx \frac{x(t_i) - x(t_{i-1})}{t_i - t_{i-1}}.
\end{equation}}

In the control algorithm, both the desired and measured strain-gauge voltages are represented by their truncated Fourier series
\begin{equation}
x(t) \approx \frac{A_0}{2} + \sum_{k = 1}^{N} \left( A_k \cos (k \omega t) + B_k \sin (k \omega t) \right),
\end{equation}
\begin{equation}
x^\ast(t) \approx \frac{A_0^\ast}{2} + \sum_{k = 1}^{N} \left( A_k^\ast \cos (k \omega t) + B_k^\ast \sin (k \omega t) \right).
\end{equation}
In our test, the first seven harmonics were retained ($N = 7$).

As a result, the total forcing $F(t)$ also can be expressed in a similar form
\begin{equation}
F(t) = \frac{A_{F0}}{2} + \sum_{k = 1}^{N} \left( A_{Fk} \cos (k \omega t) + B_{Fk} \sin (k \omega t) \right),
\end{equation}
where the coefficients are given by
\begin{subequations}
\begin{align}
  A_{F1} & = A+k_{\rm p}(A_1^\ast-A_1)-k_{\rm d} \omega (B_1^\ast-B_1), \label{eq:A1}\\
  B_{F1} & = B+k_{\rm p}(B_1^\ast-B_1)+k_{\rm d} \omega (A_1^\ast-A_1), \label{eq:B1}\\
  A_{Fk} & = k_{\rm p}(A_k^\ast-A_k)-k_{\rm d} k \omega (B_k^\ast-B_k) \; {\rm for} \nonumber \\ \; & k = 0,2,3,\dots,N, \\
  B_{Fk} & = k_{\rm p}(B_k^\ast-B_k)+k_{\rm d} k \omega (A_k^\ast-A_k)  \; {\rm for} \nonumber \\ \; & k = 2,3,\dots,N.
\end{align}
\end{subequations}

Studying the expressions \eqref{eq:A1} and \eqref{eq:B1} reveals that the fundamental harmonic component of the total forcing $\Phi = \sqrt{A_{F1}^2+B_{F1}^2}$ is not fully determined by the open-loop forcing coefficients $A$ and $B$, as it also depends on the control target $x^\ast(t)$ and the response $x(t)$. Moreover, the control introduces higher-harmonic components to the total forcing that have to be eliminated to capture the response of the open-loop system.

In our study, control-based continuation is used to generate the family of steady-state solutions of the system across a range of the forcing amplitudes while keeping the forcing frequency constant. Since, with keeping the forcing frequency $\omega$ constant, a unique forcing amplitude corresponds to every the vibration amplitude, it is possible to trace the whole branch of solutions by a sweep in the target fundamental harmonic amplitude $B_1^\ast$, with keeping $A_1^\ast = 0$ to fix the phase of the response. This means that there is a linear relationship between the continuation parameter and the forcing $F(t)$, which enables us to use a simplified version of the `full' continuation algorithm (see  \cite{Sieber2008} for example). Note that to retrieve the frequency response at a constant forcing amplitude the full algorithm would be required which may be less robust to noise than the simplified version. This issue is addressed by Schilder et al. in \cite{SCHILDER2015} with techniques developed specifically to cope with noise. {Another alternative is presented in \cite{RensonGaussian} where the full continuation algorithm is used on a local Gaussian process regression model.} Nevertheless, the frequency dependence could be equally well characterised by tracing the response curve at several frequencies as indicated by Fig. \ref{fig:sw_cbc_surf}, where response surfaces above the forcing frequency -- forcing amplitude planes, obtained by the open-loop and control-based approach, are compared. The simplified control-based continuation algorithm is briefly described below --- a full description is given in \cite{Renson2019}.  { It is also worth mentioning that, for parameter identification, it may be acceptable for the controller to be invasive, allowing it to overcome the issues around the Newton iterations. Our study however, along with \cite{SCHILDER2015} and \cite{Renson2019}, specifically focuses on recovering the bifurcation diagram, while parameter identification is used as a tool to quantify the effectiveness of the open-loop and control-based approaches.}

Let us assume that the experiment is running at a steady-state given by the target coefficients $(A_{F1}^j, B_{F1}^j, A_0^{\ast j},  B_1^{\ast j}, A_k^{\ast j}, B_k^{\ast j})$, $k = 2, \dots N$, with $A_1^{\ast j} = 0$. Then, to find the next point in the solution branch, the fundamental harmonic coefficient of the control target is increased { by the desired increment $\Delta$ as} $B_1^{\ast j+1} := B_1^{\ast j}+ \Delta$. After waiting for the control to reach steady-state, if necessary, we apply fixed point iteration to correct the higher-harmonic coefficients of the control-target until coefficients corresponding to the higher harmonics of the forcing ($A_{Fk}, \, B_{Fk}, \, k = 0,2,3, \dots N$) are below a pre-defined tolerance. Once the higher harmonics in the forcing are eliminated, the actual state, given by $(A_0^{\ast j+1},  B_1^{\ast j+1}, A_k^{\ast j+1}, B_k^{\ast j+1})$, $k = 2, \dots N$, is accepted as the steady-state response of the open loop system corresponding to the forcing amplitude
{
\begin{equation}
\Phi^{j+1} = \sqrt{\left( A_{F1}^{j+1} \right) ^2 +\left( B_{F1}^{j+1} \right)^2 }.
\end{equation}}

Provided that appropriately chosen control gains are used, this algorithm ensures a stable, non-invasive control, which traces the solution branch sweeping across the vibration amplitudes. A possible alternative could be to use a secant prediction to provide an initial guess for the algorithm in the direction obtained from the previous two points on the branch. This method may result in the algorithm reaching a fixed point in fewer iteration steps in an experiment with low noise; however, the amplitude sweep is more robust against noise since with this assumption, noise cannot affect the direction along the branch in which the next branch point is predicted which is an effect that can hinder progress along a branch. A further advantage is that the correction of the solution is carried out in a derivative-free way. Thus, it requires less evaluation at each iteration step, leading to faster convergence.

{ Note that while the control feedback given in \eqref{eq:control} is in real time in the experiment, there is no such requirement for the continuation algorithm and the setting of new control targets for the controller. Therefore, these tasks were carried out by a PC, which was also used to process the acquired data, rather than the real-time controller.}

\subsection{Model of the nonlinear oscillator}

We use the experimentally acquired bifurcation diagrams to identify the parameters of our model for the experimental rig, the one-degree-of-freedom nonlinear oscillator (see Eq. \eqref{eq:1dof_eom}). We consider a linearly damped, Duffing-like oscillator with the equation of motion
\begin{equation}\label{eq:7o_mod}
\ddot x + b \dot x + \omega_{\rm n}^2 x + \mu x^3 + \nu x^5 + \rho x^7 = \delta_{\rm st} \omega_{\rm n}^2 \cos (\omega t),
\end{equation}
where $\omega_{\rm n}$ is the linear natural angular frequency, $\delta_{\rm st}$ is the equivalent static deflection for the forcing amplitude $\delta_{\rm st} \omega_{\rm n}^2$ { (the resulting deflection of $\omega = 0$)}, while the damping is given by the parameter $b$, whereas $\mu$, $\nu$ and $\rho$ characterise the nonlinearities in the system. { Although this model does not describe the physics-based modelling of the restoring force, our investigation indicates that considering the odd nonlinear terms up to seventh order in the Duffing-type model is satisfactory to characterise the response of the experimental rig to periodic excitation. The detailed experimental characterisation of the electromagnetic forces by means of a physics-based model of the device was performed in \cite{Cammarano2011}.}

The fundamental harmonic component of the steady-state system response can be given as $X \cos (\omega t + \vartheta)$ with amplitude $X$ and phase angle $\vartheta$. Using the method of multiple scales \cite{Nayfeh}, an analytical approximate solution can be obtained for the fundamental harmonic component { if the nonlinearity is weak}. Based on this, for a given vibration amplitude $X$, the phase angle and the static deflection can be given as
\begin{equation}\label{eq:real_phaseangle}
\vartheta = \arctan \left(\frac{\tilde b \zeta}{(\zeta^2-1) - \frac{35}{64} X^6 \tilde \rho - \frac{5}{8} X^4 \tilde \nu - \frac{3}{4} X^2 \tilde \mu} \right),
\end{equation}
\begin{equation}\label{eq:real_amplitude}
\delta_{\rm st} = \left\rvert \frac{\frac{35}{64} X^7 \tilde \rho + \frac{5}{8} X^5 \tilde \nu + \frac{3}{4} X^3 \tilde \mu - X (\zeta^2-1)}{\cos(\vartheta)}\right\rvert,
\end{equation}
with $\zeta = \omega / \omega_{\rm n}$ $\tilde b := b/\omega_{\rm n}$, $\tilde \mu := \mu/\omega_{\rm n}^2$, $\tilde \nu := \nu/\omega_{\rm n}^2$, $\tilde \rho := \rho/\omega_{\rm n}^2$. The derivation of these formulae is given in the appendix.

Substituting \eqref{eq:real_phaseangle} into Eq. \eqref{eq:real_amplitude}, we obtain the static deflection by means of the system and forcing parameters as well as the amplitude of the fundamental harmonic component of the steady-state response
\begin{equation}\label{eq:static_defl}
  \delta_{\rm st} = \delta_{\rm st}(X, \tilde \mu, \tilde \nu, \tilde \rho, \tilde b, \zeta).
\end{equation}

\subsection{Numerical collocation}

\begin{figure*}[h!]
    \centering
    \includegraphics{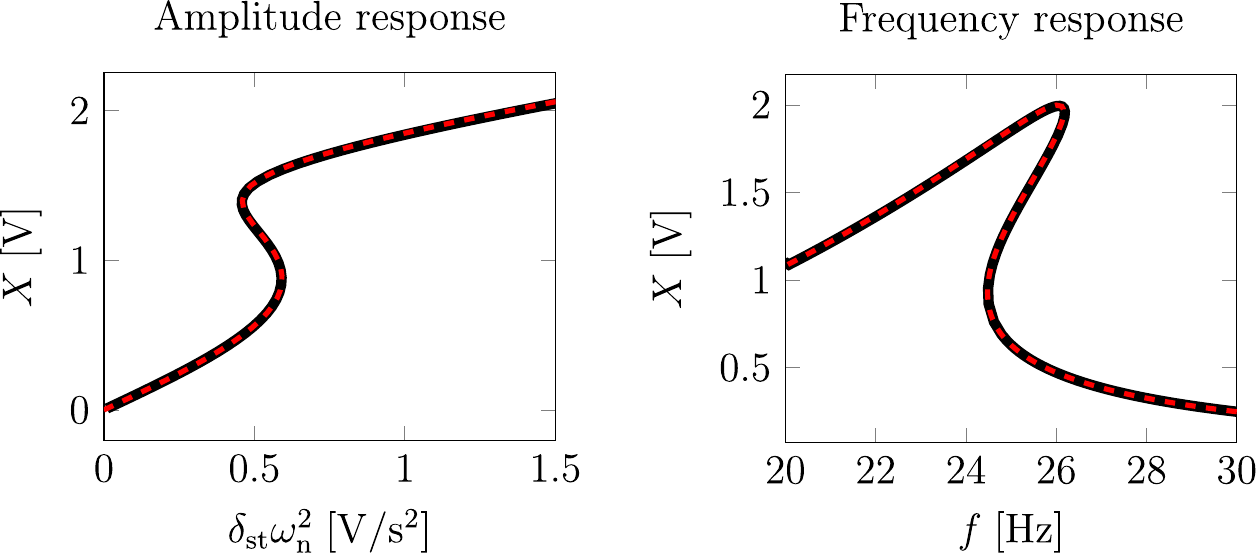}
    \caption{Comparison of the analytical approximate, and numerical solutions of Eq. \eqref{eq:7o_mod}. In the left panel, the forcing frequency is kept constant $f = 24 \, {\rm Hz}$ while the excitation amplitude is varied. In the right panel, the frequency response is shown for a constant forcing amplitude of $\delta_{\rm st} \omega_{\rm n}^2 = 0.7 \, {\rm V/s^2}$.}
    \label{fig:cbc_num}
\end{figure*}

It has to be noted that the solution presented above is only accurate for `weakly nonlinear' systems where the nonlinear terms do not dominate over the underlying linear system. To check the accuracy of the approximate solution we carried out the numerical continuation of the periodic solutions in \eqref{eq:7o_mod}. The results are compared in Fig. \ref{fig:cbc_num} for the parameters $\mu = 1.499$, $\nu = -0.3921$, $\rho = 0.0422$, $b = 0.3159$ and $f_{\rm n} = 19.95 \, {\rm Hz}$. Both frequency and amplitude variation is checked. The results indicate that the analytical approximation provides very accurate results in the parameter-range of our interest.

\section{Parameter identification}

We use the expression in Eq. \eqref{eq:static_defl} to identify the system parameters by fitting it to measurement results, based on least square errors. However, the amplitude $\delta_{\rm st} \omega_{\rm n}^2$ of the forcing provided by the shaker is not measured directly. Instead, we capture the acceleration $a_{\rm base}$ of the base of the oscillator. We assume that the base acceleration and forcing are proportional: $\delta_{\rm st} = c_A A_{\rm base}$, where $A_{\rm base}$ is the fundamental harmonic amplitude of the base acceleration. This leads to an expression for the base acceleration amplitude that can be directly used for parameter identification
\begin{equation}\label{eq:volt_amp}
  A_{\rm base} = \frac{1}{c_A} \delta_{\rm st}(X, \zeta, \tilde \mu, \tilde \nu, \tilde \rho, \tilde b).
\end{equation}
Note that the base acceleration is not controlled directly in the experiment as the real time controller sets the shaker voltage instead. In Eq. \eqref{eq:volt_amp}, the frequency ratio $\zeta$ is set by the user (provided the natural angular frequency $\omega_{\rm n}$ is a priori determined) whereas the response amplitude $X$ is measured. Consequently, we aim to identify the parameters $\tilde b$, $\tilde \mu$, $\tilde \nu$, $\tilde \rho$ and $c_A$ such that they minimise the least square error based on Eq. \eqref{eq:volt_amp}. Assuming that $m$ measurement points, given by $({A_{\rm base}}_i,\zeta_i, X_i)$, are acquired, the function to be minimised can be expressed as
\begin{equation}\label{eq:objective}
  R = \sum_{i=1}^{m} \left( {A_{\rm base}}_i - \frac{1}{c_A} \delta_{\rm st} (X_i, \zeta_i, \tilde \mu, \tilde \nu, \tilde \rho, \tilde b ) \right)^2.
\end{equation}
The least square fit is obtained by using the gradient-free principal-axis method of the NLopt package in Julia \cite{NLOPT}. Once the least square fit is performed, the original model parameters $b = \tilde b \omega_{\rm n}$, $\mu = \tilde \mu \omega_{\rm n}^2$, $\nu = \tilde \nu \omega_{\rm n}^2$, $\rho = \tilde \rho \omega_{\rm n}^2$ can be recovered.

{ Note that due to the nonlinear nature of the fitting problem, there is no strong evidence that the identified parameters correspond to the global minimum of the objective function in \eqref{eq:objective}. Nevertheless, the visual assessment of the response of the filled model suggest that this method provides a good fit to the measurement data. It is also worth mentioning that one may include the formula in \eqref{eq:real_phaseangle} for the phase-angle $\vartheta$ in the parameter-identification. Nevertheless, since the phase-lag is changing monotonously with the vibration amplitude ($\vartheta \approx -\pi$ for small amplitudes and $\vartheta \approx 0$ at large amplitudes) along the S-shaped curve, eliminating it from the objective function does not change the main tendencies we observe in this study.}
\section{Effect of noise on parameter identification}

\subsection{Underlying linear response}

\begin{figure}
\begin{center}
\includegraphics{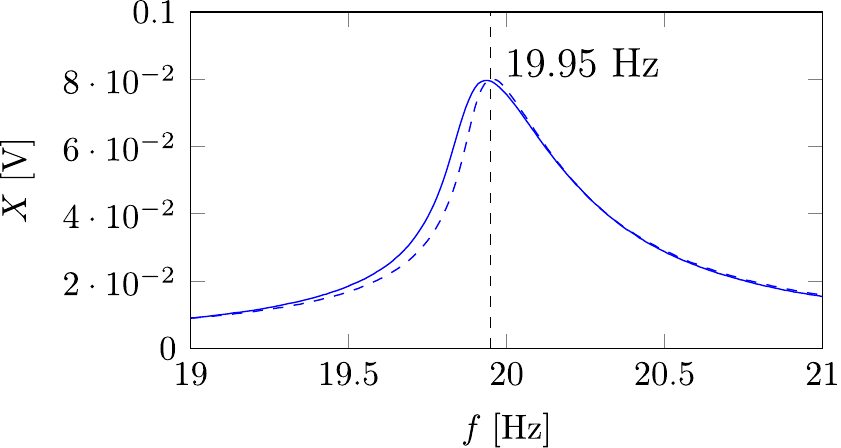}
\caption{Linear frequency response: up-sweep (blue continuous), down-sweep (blue, dashed)}\label{fig:lin_resp}
\end{center}
\end{figure}

The measurement procedure is as follows. Firstly, the linear frequency response of the system is extracted by performing an open-loop frequency sweep with a low shaker voltage amplitude ($U_{\rm shaker} = 0.01 \; \rm V$) where the effect of the nonlinearities is marginal. We performed an up and a down-sweep between 19-21 Hz with an increment of $0.1 \; \rm{Hz}$ (see Fig. \ref{fig:lin_resp}). One can observe a slight discrepancy between the two curves which can be explained by the effect of temperature change during the measurement. Estimating the damping with the 3 dB bandwidth method indicated that, at this level of excitation, the Lehr's damping ratio of the system is below 0.01. This means that the difference between the measurable `damped' and the undamped natural frequencies is below $0.01 \%$. Therefore, we accepted the locus of the maximum of a response curve as the undamped natural frequency. Since the up- and down-sweeps provided different results, we calculate with the average of the two peaks $f_{\rm n} = 19.95 \; \rm{Hz}$ as the natural frequency of the system.

\subsection{Robustness of solution tracing methods against process noise}

\begin{figure*}[]
  \begin{centering}
    \includegraphics{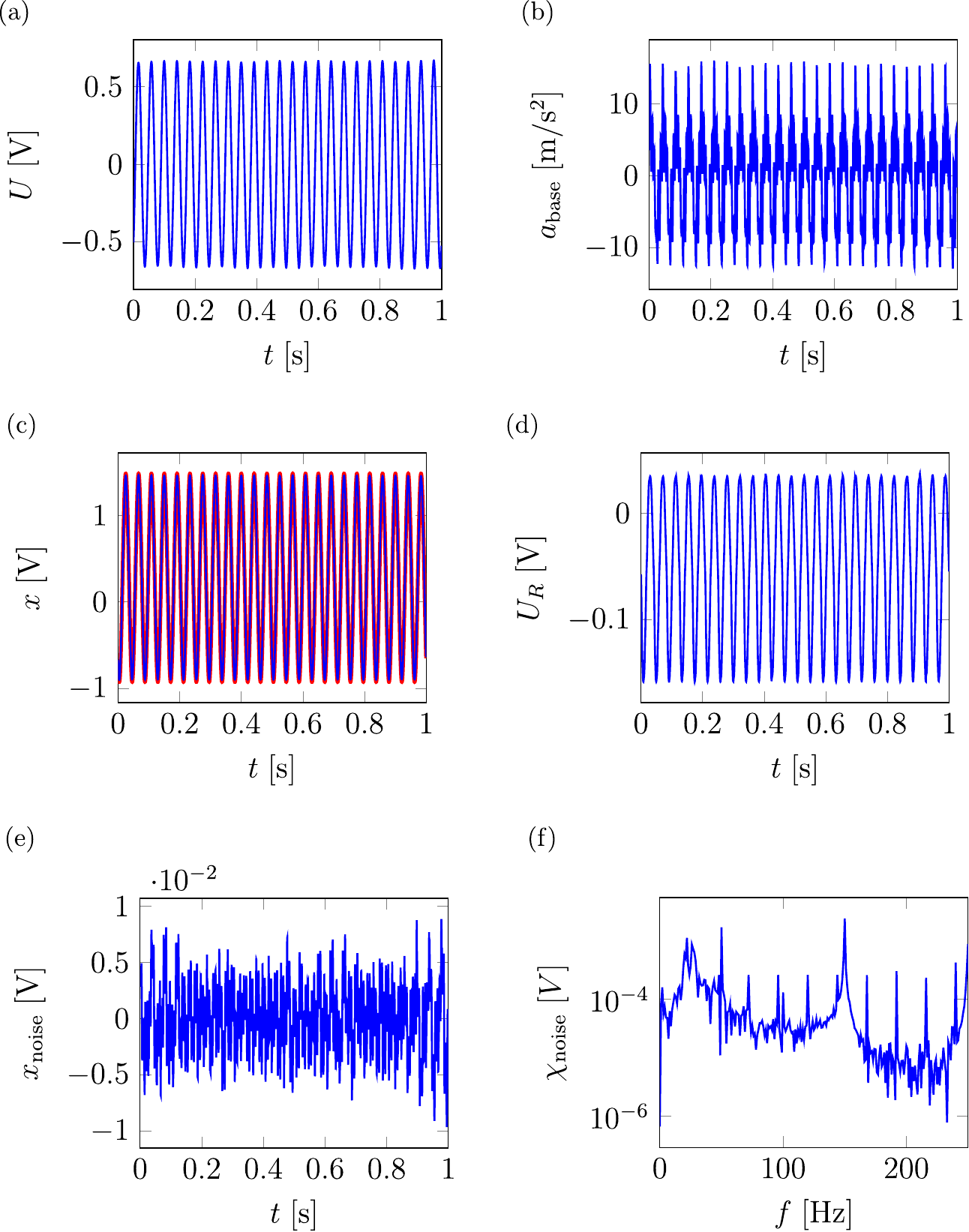}
    \caption{Acquired data and outputs generated by the real time controller at low noise-level. (a) shaker voltage $U$, (b) base acceleration $a_{\rm base}$, (c) target (red) and actual (blue) strain-gauge voltage $x$, (d) resistor voltage $U_{\rm R}$, (e) noise time profile $x_{\rm noise}$, (f) noise FFT $\chi_{\rm noise}$.}\label{fig:signals_no_noise}
  \end{centering}
\end{figure*}

\begin{figure*}[]
  \begin{centering}
	\includegraphics{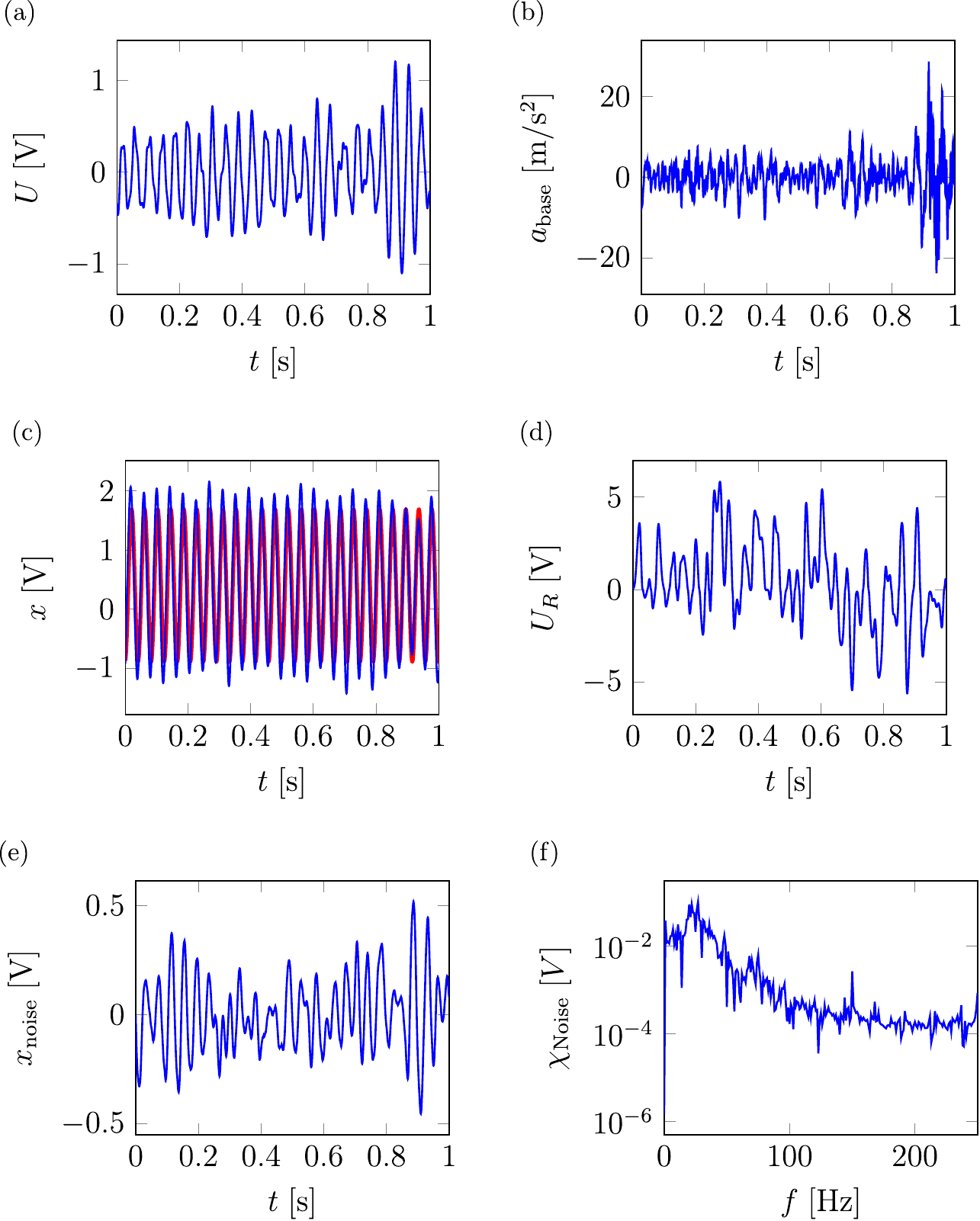}
    \caption{Acquired data and outputs generated by the real time controller at high noise-level (N10). (a) shaker voltage $U$, (b) base acceleration $a_{\rm base}$, (c) target (red) and actual (blue) strain-gauge voltage $x$, (d) resistor voltage $U_{\rm R}$, (e) noise time profile $x_{\rm noise}$, (f) noise FFT $\chi_{\rm noise}$.}\label{fig:signals_with_noise}
  \end{centering}
\end{figure*}

The core of our study was to assess the performance of control-based continuation against open-loop parameter-sweeps in scenarios where the experiment was polluted with different levels of noise. We used the real-time controller to generate a random demand for the current flowing through the electromagnetic coil by generating a series of pseudo-random numbers and feeding the signal through a low-pass Butterworth filter resulting in a band-limited noise with a cutoff frequency of 50 Hz \cite{DBrtc}. Time profiles of the input and output signals, without and with additional noise, are shown in Fig. \ref{fig:signals_no_noise} and \ref{fig:signals_with_noise} respectively. The noise $x_{\rm noise}$ in these diagrams is obtained by subtracting the harmonic components of the strain-gauge voltage $x$, identified by the real-time controller.

In the case where no noise is added through the coil, see Fig. \ref{fig:signals_no_noise}, other uncontrolled noise sources provide a low noise-level of below 1\% of the vibration amplitude. The notable peaks in the noise FFT belong to the power supply frequency and its upper harmonics (50, 100, 150 Hz) while one can also observe smaller peaks ($n \times 24 \; {\rm [Hz]}$) related to numerical errors in the calculation of the harmonic coefficients. These components become negligible when significant noise is added to the system through the coil, resulting in more evenly distributed frequency components as shown in Fig. \ref{fig:signals_with_noise}.

While the control is clearly affected by the added noise, the algorithm is still capable to stabilise the system around the periodic solutions of the underlying noise-free system even when the noise amplitude is about 30\% of the vibration amplitude.

Also note how the presence of noise affects the current in the coil (measured through the voltage on a resistor): if there is no additional noise in the system, the coil current appears to be largely harmonic, as its main source is electromagnetic induction due to the tip mass passing the iron core. This component practically disappears in the additional noise.

Using control-based continuation in a noise-polluted environment is more challenging as error tolerances are more difficult to meet. Moreover, if a secant predictor is used, as is common with the pseudo-arclength method, the algorithm may fail resulting in repeated, or to the contrary, no coverage of some parts of the solution manifold.

In our case though, it proved to be sufficient to adjust the simplified control-based continuation algorithm by averaging the coefficients of 10 periods and allowing a higher error tolerance in the noisy measurements. In some cases though, even increased tolerances were difficult to meet; therefore, at higher noise levels we omitted the step of the algorithm that corrects the higher harmonics and accepting the solution we obtained by simply changing the fundamental harmonic component of the control target. It is worth mentioning that this is always the case in the open-loop measurements when steady-state solutions are accepted without any correction step.

We also used the fact that, due to the choice of continuation in amplitude (rather than frequency), a unique solution exists for every response amplitude. Thus, performing a sweep in the targeted vibration amplitudes ensured that the continuation of the solution branch did not stall and the whole branch was covered.

\begin{figure}[]
    \centering
    	\includegraphics{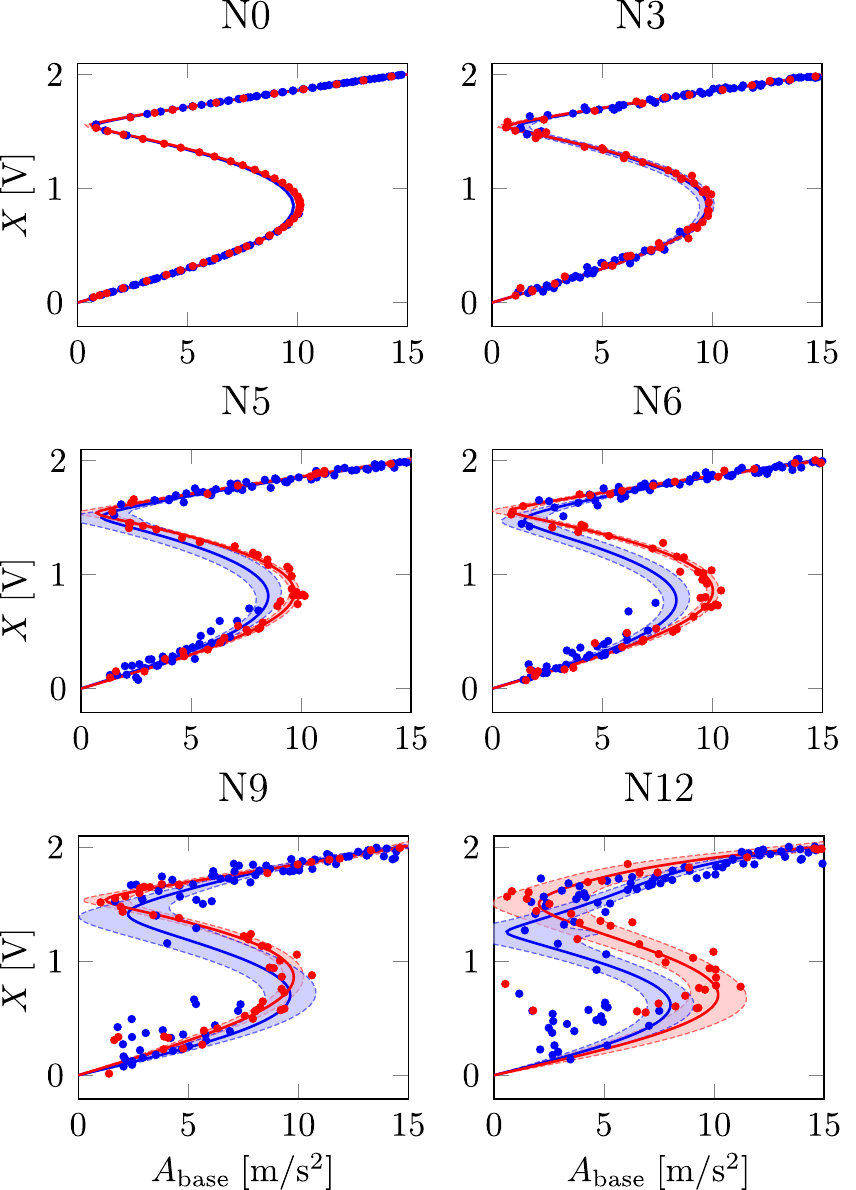}
\caption{Amplitude response of the nonlinear oscillator at 24 Hz, at different noise levels. The markers show the measurement data while the response of the fitted models is indicated by the thick continuous lines. {The shaded areas and the dashed boundaries indicate the $95\%$ approximate confidence intervals for the identified response curves.} The blue markers and curves belong to parameter-sweeps while the red ones to CBC.}
\label{fig:cbc_sw_fits}
\end{figure}

\begin{figure*}[]
    \centering
    \includegraphics{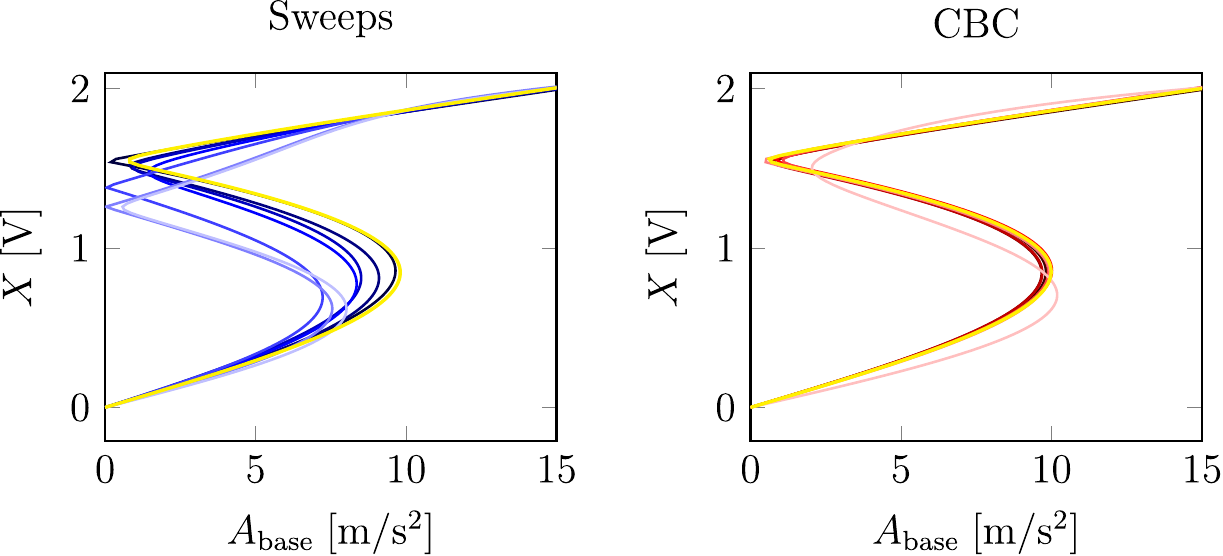}
    \caption{Robustness against noise: The amplitude responses fitted to amplitude sweeps (left panel) and control based continuation (right panel) with base acceleration. The curves belong to noise-levels N2, N4, N5, N6, N8, N10 and N12 with higher noise indicated by lighter shading. The response-curves belonging to the cases without additional noise (N0) are highlighted in yellow.}
    \label{fig:noise_robustness}
\end{figure*}

In Fig. \ref{fig:cbc_sw_fits}, we compared control-based continuation and open-loop forcing amplitude sweeps under different levels of noise load. The labels beginning with N indicate the magnitude of noise we added through the coil, e.g. N3 means that the noise magnitude is three times as the reference N1 while N0 corresponds to the case when the experiment was not polluted with noise through the coil.

{In the panels N0, N3 and N5 of Fig. \ref{fig:cbc_sw_fits}, we show the measurements where we had a maximum tolerance for the higher harmonics in the control-based continuation algorithm whereas in the measurements shown in the panels N6, N9 and N12}, we accepted the result without introducing corrections to the higher harmonic components. Comparing the cases N5 and N6 which were performed with a similar amount of noise but one with, and the other without higher harmonic control, we can see that we did not introduce significant further error in the experiment by not setting a maximum tolerance at higher noise levels. { Note that, even in the case without additional noise, we recorded a few measurement points close to the unstable part of the S-shaped curves with the open-loop method, even though this would not be possible without control. This apparent contradiction can be explained by the fact that we controlled the input voltage of the shaker in the experiment instead of the base acceleration. The reason for using the measured acceleration as an input parameter is that we found that treating the actuator as part of our system would significantly change its dynamics. As a result, instead of seeing a sudden drop in the vibration amplitude, as one would expect in the open-loop case, we found that the transient behaviour, at least initially, kept the system close to the weakly unstable region of the branch, and our algorithm identified the slowly diverging solutions as `periodic' due to the averaging, which was being used to partly mitigate the effects of noise.}

\begin{table*}
\begin{center}
{
\begin{tabular}{|c|cc|cc|cc|cc|cc|}\hline
  \multicolumn{11}{|c|}{Open-loop measurements}\\
  \hline
  \multirow{2}{*}{Noise} & \multicolumn{2}{c|}{$\tilde \mu$} & \multicolumn{2}{c|}{$\tilde \nu$} & \multicolumn{2}{c|}{$\tilde \rho$} & \multicolumn{2}{c|}{$\tilde b$} & \multicolumn{2}{c|}{$c_U$}\\
  & IDV & ESTD & IDV & ESTD & IDV & ESTD & IDV & ESTD & IDV & ESTD
  \\\hline
  N0 & 0.2997 & 0.0055 & -0.0230 & 0.0037 & -0.00055 & 0.00056 & 0.01308 & 0.00201 & 0.02669 & 0.00014 \\
  N2 & 0.3039 & 0.0122 & -0.0268 & 0.0082 & -8.278e-7 & 0.00126 & 1.461e-7 & 249.94 & 0.02590 & 0.00031\\
  N3 & 0.3047 & 0.0148 & -0.0273 & 0.0100 & -6.778e-5 & 0.00153 & 0.01589 & 0.00377 & 0.02766 & 0.00040\\
  N4 & 0.3495 & 0.0295 & -0.0556 & 0.0192 & 0.00424 & 0.00290 & 0.01146 & 0.01001 & 0.02580 & 0.00073\\
  N5 & 0.3467 & 0.0292 & -0.0526 & 0.0213 & 0.00402 & 0.00320 & 0.01403 & 0.01241 & 0.02766 & 0.00081\\
  N6 & 0.3990 & 0.0428 & -0.0867 & 0.0195 & 0.00886 & 0.00300 & 0.02175 & 0.00793 & 0.02662 & 0.00094\\
  N8 & 0.5165 & 0.0448 & -0.1636 & 0.0293 & 0.02099 & 0.00446 & 1.819e-10 & 303425.5 & 0.02717 & 0.00163\\
  N9 & 0.5068 & 0.0438 & -0.1602 & 0.0289 & 0.01934 & 0.00446 & 0.02746 & 0.01310 & 0.02075 & 0.00134\\
  N10 & 0.6366 & 0.0533 & -0.2413 & 0.0321 & 0.03196 & 0.00463 & 4.928e-9 & 499541.7 & 0.02340 & 0.00169\\
  N12 & 0.6448 & 0.0379 & -0.2488 & 0.0269 & 0.03319 & 0.00441 & 0.00830 & 0.03163 & 0.02191 & 0.00153\\\hline
   \multicolumn{11}{|c|}{Control-based continuation}\\
  \hline
 \multirow{2}{*}{Noise} & \multicolumn{2}{c|}{$\tilde \mu$} & \multicolumn{2}{c|}{$\tilde \nu$} & \multicolumn{2}{c|}{$\tilde \rho$} & \multicolumn{2}{c|}{$\tilde b$} & \multicolumn{2}{c|}{$c_U$}\\
  & IDV & ESTD & IDV & ESTD & IDV & ESTD & IDV & ESTD & IDV & ESTD
  \\\hline
  N0$^\ast$ & 0.2999 & 0.0035 & -0.0258 & 0.0025 & -0.00025 & 0.00038 & 0.00798 & 0.00184 & 0.02496 &  0.00013\\
  N2$^\ast$ & 0.2989 & 0.0064 & -0.0233 & 0.0047 & -0.00067 & 0.00075 & 0.01004 & 0.00225 & 0.02538 & 0.00024\\
  N3$^\ast$ & 0.2929 & 0.0106 & -0.0186 & 0.0078 & -0.00158 & 0.00129 & 0.01018 & 0.00340 & 0.02538 & 0.00037\\
  N4$^\ast$ & 0.3132 & 0.0132 & -0.0344 & 0.0098 & 0.00114 & 0.00159 & 0.01217 & 0.00554 & 0.02538 & 0.00047\\
  N5$^\ast$ & 0.3055 & 0.0168 & -0.0270 & 0.0132 & -0.00031 & 0.00231 & 0.00970 & 0.00648 & 0.02549 & 0.00054\\
  N6$^\dagger$ & 0.2917 & 0.0187 & -0.0206 & 0.0135 & -0.00102 & 0.00215 & 0.01124 & 0.00608 & 0.02516 & 0.00065\\
  N8$^\dagger$ & 0.3032 & 0.0237 & -0.0276 & 0.0165 & -0.00010 & 0.00252 & 0.01425 & 0.00540 & 0.02512 & 0.00084\\
  N9$^\dagger$ & 0.2755 & 0.0316 & -0.0017 & 0.0237 & -0.00493 & 0.00392 & 0.01760 & 0.00734 & 0.02605 & 0.00104\\
  N10$^\dagger$ & 0.3160 & 0.0292 & -0.0365 & 0.0210 & 0.00130 & 0.00332 & 0.00578 & 0.01711 & 0.02460 & 0.00107\\
  N12$^\dagger$ & 0.5214 & 0.0585 & -0.1889 & 0.0444 & 0.02604 &  0.00753 & 0.02189 & 0.01115 & 0.01947 & 0.00186\\\hline
   \end{tabular}}
  \caption{{The identified model parameter values (IDV) and estimated standard deviations (ESTD) with open-loop amplitude-sweeps and control-based continuation at different noise levels: the coefficients $\tilde \mu, \; \tilde \nu \; \tilde \rho$ of the $3^{\rm rd}$, $5^{\rm th}$ and $7^{\rm th}$ order terms, the linear viscous damping $\tilde b$ and the scaling factor $c_{\rm A}$ between the forcing and base acceleration. In case of the control-based continuation data, the symbols $\ast$ and $\dagger$ refer to measurements with and without a maximum tolerance for the higher harmonic components of the forcing, respectively.}}\label{tab:sw}
\end{center}
\end{table*}

It can be observed that if the experiment is running with a low amount of noise (see N3) both methods, the standard open-loop sweep and control-based continuation provide a response close to the case with no additional noise (N0), which we take to be the baseline `correct' solution. However, at medium noise levels (see N5, N6 and N9) the jumps in the open-loop data between the low and high amplitude parts of the solution branch occur significantly earlier than the folds in the branch of steady-state solutions of the underlying noise-free system. Meanwhile, control-based continuation is still able to trace the entire family of solutions and to retain the response curve of the low noise system with reasonable accuracy. This is true even at high noise level (see N12) to some extent --- even though control-based continuation begins to struggle in tracking the branch, it still provides some information about the unstable solutions while the bistable region cannot be recovered by parameter sweeps.

This can be well demonstrated by Fig. \ref{fig:noise_robustness} where we compare the response curves obtained by fitting the parameters of the seventh-order Duffing-like oscillator, given by Eq. \eqref{eq:7o_mod}, to measurement results acquired at different noise levels.

{ The robustness of the parameter identification is assessed by linearising the model in \eqref{eq:volt_amp} around the identified parameters

\begin{equation}
    \begin{split}
  A_{\rm base} & \approx \frac{1}{c_A^\ast} \delta_{\rm st}(X, \zeta, \tilde \mu^\ast, \tilde \nu^\ast, \tilde \rho^\ast, \tilde b^\ast) \\
  & + \mathbf{grad}_{\mathbf p} \left. \left( \frac{1}{c_A} \delta_{\rm st}(X, \zeta, \tilde \mu, \tilde \nu, \tilde \rho, \tilde b) \right) \right\rvert_{\mathbf{p} = \mathbf{p}^\ast} \mathbf h, \\
\end{split}
\end{equation}
where the parameter vector reads $\mathbf{p} = (\tilde \mu, \tilde \nu, \tilde \rho, \tilde b, c_A)^{\rm T}$, $\mathbf{p}^\ast$ is the minimiser of the objective fuction \eqref{eq:objective}, while $\mathbf{h} = \mathbf{p} - \mathbf{p}^\ast$.

Then, the standard deviations of the identified parameters are estimated by means of the asymptotic distribution \cite{nonlinreg}. The identified parameter values and their estimated standard deviation is listed in Table \ref{tab:sw} for both sets of measurements, amplitude-sweeps and control-based continuation, respectively. Note that, in case of the damping parameter $\hat b$, our model is singular at $\hat b=0$. Thus, the estimation procedure provided unrealistic values for the standard deviation for values of $\hat b \approx 0$ (see the noise levels N2, N8 and N10 in the open-loop case).  Moreover, we used the linearised model to draw the $95\%$ confidence intervals around the identified S-shaped response curves in Fig. \ref{fig:cbc_sw_fits} as a visual indicator of the robustness of parameter fitting.

From these results, it is clear that control-based continuation is able to preserve the original response curve better than frequency sweeps. Furthermore, at low and medium noise levels (N0-N9) the S-shaped curves, obtained from the control-based continuation data, are not just better representing the deterministic dynamics but they fit more robustly to the measurement data than in case of the open-loop measurements. This changes at high noise-level (case N12), where the open loop data results in lower standard deviation in the identified parameters and the confidence intervals of the identified response curves are comparable for the open-loop and the control-based measurements. This feature of the control-based continuation data can also be seen as an advantage, as this indicates that, at this noise level, there is a significant loss of information about the underlying system even with control-based continuation. In the meantime, fitting to the open-loop parameter sweeps at this noise level results in a model that seems to be more robust according to these measures even though it is clearly less accurate when compared to the response with no additional noise.}

The identified model parameters are listed in Table \ref{tab:sw} for both set of measurements, the amplitude-sweeps and control-based continuation, respectively.

%\begin{table}
%\begin{center}
%  \begin{tabular}{cccccc}\hline
%    $c_{\rm noise}$ & $\hat \mu$ & $\hat \nu$ & $\hat \rho$ & $\hat b$ & $c_U$ \\\hline
%    N0$^\ast$ & 0.29995 & -0.0257753 & -0.000252701 & 0.0079761 & 0.0249586\\
%    N2$^\ast$ & 0.298945 & -0.0233512 & -0.000669121 & 0.010041 & 0.0253846\\
%    N3$^\ast$ & 0.292901 & -0.0185704 & -0.00157755 & 0.0101766 & 0.0253763\\
%    N4$^\ast$ & 0.313156 & -0.0344051 & 0.00113597 & 0.0121747 & 0.0253811\\
%    N5$^\ast$ & 0.305511 & -0.0269727 & -0.000312632 & 0.00969982 & 0.0254893\\
%    N6$^\dagger$ & 0.291671 & -0.0206321 & -0.00101879 & 0.0112381 & 0.0251575\\
%    N8$^\dagger$ & 0.303191 & -0.027599 & -0.000102184 & 0.0142532 & 0.0251187\\
%    N9$^\dagger$ & 0.275549 & -0.00172186 & -0.0049329 & 0.0175988 & 0.0260507\\
%    N10$^\dagger$ & 0.315962 & -0.0365078 & 0.00129476 & -0.00577997 & 0.0246026\\
%    N12$^\dagger$ & 0.521444 & -0.188949 & 0.0260422 & 0.0218853 & 0.0194733\\\hline
%  \end{tabular}
%  \caption{The identified model parameters with control-based continuation in forcing amplitude at different noise levels. The symbols and $\ast$ and $\dagger$ refer to measurements with and without a maximum tolerance for the higher harmonic components of the forcing, respectively.}\label{tab:cbc}
%\end{center}
%\end{table}

\section{Conclusions}

By investigating the steady-state response of a forced nonlinear oscillator under different levels of process noise, the robustness of control-based continuation was assessed by comparing it to open-loop measurements.

Based on the acquired data, we carried out parameter identification of a seventh-order Duffing-like oscillator to quantify and characterise the performance of the two solution tracking methods. We demonstrated that the ability of control-based continuation to capture both stable and unstable periodic solutions, and the fact that we have feedback control on the response, result in a more robust coverage of the solution branch than in case of parameter-sweeps. Furthermore, while an open-loop parameter sweep can perform as well as control-based continuation at low noise levels, it has a tendency to miss parts of the stable solutions in bistable parameter domains as the additional noise can cause the system to jump between the domains of attraction of the co-existing stable solutions. This can result in poorly identified parameters in the model, while fitting the model to data obtained by control-based continuation yields a response-curve that is closer to the response of the underlying noise-free system.

As control-based continuation is affected by process noise, albeit to a lesser degree, it is still an open question to what extent can one expect recover the response of the underlying deterministic system in a heavily noise-contaminated measurement. In general, control-based continuation is capable to extract more information from experiments than open-loop measurements. Thus, with the help of control-based continuation it may be possible to capture finer, otherwise undetectable, details about the dynamics of physical systems.

\section*{Data statement}

All the experimental data used in this paper have been deposited into the University of Bristol
Research Data Repository and is publicly available for download.

\section*{Conflicts of interest}
The authors have no conflicts of interest to declare that are relevant to the content of this article.

\section*{Preprint}
The preprint \cite{thispreprint} has been deposited to arXiv and can be accessed at \url{https://arxiv.org/abs/2001.11008}.

\begin{acknowledgements}
This research has received funding from the {\em Next-generation test methods for nonlinear structures} (EP/P019323/1) and the {\em Digital twins for improved dynamic design} (EP/R006768/1) EPSRC grants. The support of the EPSRC is greatly acknowledged.
\end{acknowledgements}



\bibliographystyle{unsrt}
\bibliography{cbc6}   % name your BibTeX data base
% \begin{thebibliography}{}
% %
% % and use \bibitem to create references. Consult the Instructions
% % for authors for reference list style.
% % Format for Journal Reference
% Author, Article title, Journal, Volume, page numbers (year)
% Author, Book title, page numbers. Publisher, place (year)
% \end{thebibliography}

\appendix
\section*{Appendix: Analytical approximation of the response of the seventh-order forced Duffing-oscillator}

%\numberwithin{equation}{section}

To characterise our experimental rig, as a model we use a seventh-order, Duffing-like oscillator with the equation of motion given by \eqref{eq:7o_mod}:
\begin{equation}\label{eq:7o_mod_app}
\ddot x + b \dot x + \omega_{\rm n}^2 x + \mu x^3 + \nu x^5 + \rho x^7 = \delta_{\rm st} \omega_{\rm n}^2 \cos (\omega t),
\end{equation}
with a linear natural angular frequency $\omega_{\rm n}$, a static deflection $\delta_{\rm st}$, for the viscous damping $b$, and nonlinear coefficients $\mu$, $\nu$ and $\rho$.

For the steady-state response of this system, one can derive an analytical approximate solution using the method of multiple scales \cite{Nayfeh}. Accordingly, we assume that a steady-state solution of \eqref{eq:7o_mod_app} can be expanded in the following form
\begin{equation}\label{eq:ueps}
  u(t) = u_0 (t) + \varepsilon u_1 (t) + \varepsilon^2 u_2 (t) + \dots,
\end{equation}
where $u_0$ is the fundamental harmonic component of the solution, considered as a steady-state solution of the linear, undamped, homogeneous system
\begin{equation}\label{eq:Lin}
\ddot x + \omega^2 x = 0,
\end{equation}
with a natural frequency equal to the forcing frequency, while $\varepsilon$ is the perturbation parameter, which is assumed to be sufficiently small.

To simplify our formulae, we introduce dimensionless time as $\tau := \omega_{\rm n} t$. Using the transformation $\rm d/\rm d t = \omega_{\rm n} \rm d / \rm d \tau$ between the derivatives, the equation of motion \eqref{eq:7o_mod} can be expressed as
\begin{equation}
x'' + \frac{b}{\omega_{\rm n}} x' + x + \frac{\mu}{\omega_{\rm n}^2} x^3 + \frac{\nu}{\omega_{\rm n}^2} x^5 + \frac{\rho}{\omega_{\rm n}^2} x^7 = \delta_{\rm st} \cos (\omega t),
\end{equation}
where primes denote the derivation with respect to the dimensionless time $\tau$. With the aim to express this equation with the homogenous system \eqref{eq:Lin} and a nonlinear perturbation we introduce new coefficients for the nonlinear, damping and forcing terms $\hat \mu := \mu / (\varepsilon \omega_{\rm n}^2)$, $\hat \nu := \nu / (\varepsilon \omega_{\rm n}^2)$, $\hat \rho := \rho / (\varepsilon \omega_{\rm n}^2)$, $\hat b := b / (\varepsilon \omega_{\rm n})$, as well as the frequency ratio $\zeta = \omega/ \omega_{\rm n}$. This yields to
\begin{equation}
x'' + x  = \varepsilon \left(\frac{ \delta_{\rm st}}{\varepsilon} \cos (\zeta \tau) - \hat b x'- \hat \mu x^3 - \hat \nu x^5 -\hat \rho x^7 \right).
\end{equation}
Adding $\zeta^2 x$ to both sides, this equation can be re-arranged as
\begin{equation}
\begin{split}
\omega_{\rm n}^2 x'' & + \zeta^2 x  = \\
& \varepsilon \left(\dfrac{\zeta^2-1}{\varepsilon} x - \hat b x'- \hat \mu x^3 - \hat \nu x^5 -\hat \rho x^7 - \frac{\delta_{\rm st}}{\varepsilon} \cos (\zeta \tau) \right).
\end{split}
\end{equation}
Introducing the parameters $\Psi : = (\zeta^2-1)/ \varepsilon $ and $\Phi : = \hat \delta_{\rm st} / \varepsilon $ we can express the equation of motion in the perturbed form of
\begin{equation}\label{eq:ndim_eom}
\omega_{\rm n}^2 x'' + \zeta^2 x  = \varepsilon \left(\Psi x - \hat b x'- \hat \mu x^3 - \hat \nu x^5 -\hat \rho x^7 - \Phi \cos (\zeta \tau) \right).
\end{equation}
Substituting in the solution \eqref{eq:ueps} and expanding equation \eqref{eq:ndim_eom} in a power series form by means of the perturbation parameter $\varepsilon$ one obtains
\begin{subequations}
  \begin{align}
u_0'' + \zeta^2 u_0 & = 0, \label{eq:eps0}\\
u_1'' + \zeta^2 u_1 & = \Psi u_0 - \hat b u_0' - \hat \mu u_0^3 - \hat \nu u_0^5 - \hat \rho u_0^7 + \Phi \cos (\zeta \tau), \label{eq:eps1}\\
u_2'' + \zeta^2 u_2 & = \dots \, .
\end{align}
\end{subequations}
The general solution of the homogenous ODE \eqref{eq:eps0} can be given as
\begin{equation}\label{eq:u0}
  u_0(\tau) = X \cos (\zeta \tau + \vartheta),
\end{equation}
with the vibration amplitude $X$ and phase angle $\vartheta$.

Note that this solution also contributes to the `forcing' in the right-hand side of Eq. \eqref{eq:eps1}. Moreover, $u$ is a steady-state solution of \eqref{eq:7o_mod_app} which implies that the functions $u_k, k=1,2,3,\dots$ should be bounded: $\left\rvert u_k(\tau)\right\rvert \leq M$ for $\tau \in [0, \infty)$, $M \in {\rm I\!R^+}$. Thus, the amplitude $X$ and the phase angle $\vartheta$ in \eqref{eq:u0} should ensure that a non-resonant forcing is provided in Eq. \eqref{eq:eps1}.

Substituting the solution \eqref{eq:u0} in Eq. \eqref{eq:eps1} yields to
\begin{equation}\label{eq:eps1_1}
  \begin{split}
u_1'' & + \zeta^2 u_1 = \Psi X \cos (\zeta \tau + \vartheta) + \zeta \hat b X \sin (\zeta \tau + \vartheta) \\
 & - \hat \mu X^3 \cos^3 (\zeta \tau + \vartheta) - \hat \nu X^5 \cos^5 (\zeta \tau + \vartheta) \\
 & - \hat \rho X^7 \cos^7 (\zeta \tau + \vartheta) + \Phi \cos (\zeta \tau).
\end{split}
\end{equation}
The linearisation of the trigonometric expressions in the equation above provides
\begin{equation}
  \begin{split}
  u_1'' & + \zeta^2 u_1  = \\ & \left( \frac{35}{64} X^7 \hat \rho + \frac{5}{8} X^5 \hat \nu + \frac{3}{4} X^3 \hat \mu - X \Psi - \Phi \cos(\vartheta) \right) \cos(\zeta \tau + \vartheta) \\ & - \left(X \hat b \zeta + \Phi \sin (\vartheta) \right) \sin(\zeta \tau + \vartheta) + \dots,
\end{split}
\end{equation}
where we include only the base harmonic terms, since these are the terms which could potentially trigger a resonance in Eq. \eqref{eq:eps1_1}. This can be avoided only if the coefficients of $\cos(\zeta \tau + \vartheta)$ and $\sin(\zeta \tau + \vartheta)$ are zero. That is
\begin{equation}\label{eq:nonresonant_1}
\frac{35}{64} X^7 \hat \rho + \frac{5}{8} X^5 \hat \nu + \frac{3}{4} X^3 \hat \mu - X \Psi - \Phi \cos(\vartheta) = 0,
\end{equation}
\begin{equation}\label{eq:nonresonant_2}
X \hat b \zeta + \Phi \sin (\vartheta) = 0.
\end{equation}
Using equations \eqref{eq:nonresonant_1} and \eqref{eq:nonresonant_2} one can derive an expression for the phase angle and the forcing amplitude each
\begin{equation}\label{eq:phaseangle}
\vartheta = \arctan \left(\frac{\hat b \zeta}{\Psi - \frac{35}{64} X^6 \hat \rho - \frac{5}{8} X^4 \hat \nu - \frac{3}{4} X^2 \hat \mu} \right),
\end{equation}
\begin{equation}\label{eq:amplitude}
\Phi = \left\rvert \frac{\frac{35}{64} X^7 \hat \rho + \frac{5}{8} X^5 \hat \nu + \frac{3}{4} X^3 \hat \mu - X \Psi}{\cos(\vartheta)} \right\rvert.
\end{equation}
For the static deflection $\delta_{\rm st}$ these provide
\begin{equation}
\vartheta = \arctan \left(\frac{\tilde b \zeta}{(\zeta^2-1) - \frac{35}{64} X^6 \tilde \rho - \frac{5}{8} X^4 \tilde \nu - \frac{3}{4} X^2 \tilde \mu} \right),
\end{equation}
\begin{equation}
\delta_{\rm st} = \left\rvert \frac{\frac{35}{64} X^7 \tilde \rho + \frac{5}{8} X^5 \tilde \nu + \frac{3}{4} X^3 \tilde \mu - X (\zeta^2-1)}{\cos(\vartheta)}\right\rvert,
\end{equation}
with $\tilde b := \varepsilon \hat b$, $\tilde \mu := \varepsilon \hat \mu$, $\tilde \nu := \varepsilon \hat \nu$, $\tilde \rho := \varepsilon \hat \rho$.

\end{document}